\documentclass[3p,times]{elsarticle}
\usepackage{amsfonts}
\usepackage{amssymb}
\usepackage{amsmath}
\usepackage{mathrsfs}
\numberwithin{equation}{section}
\newtheorem{thm}{Theorem}[section]
\newtheorem{lem}[thm]{Lemma}
\newtheorem{prop}[thm]{Proposition}

\newdefinition{remark}[thm]{Remark}
\newdefinition{defi}[thm]{Definition}
\newproof{pf}{Proof}

\usepackage[figuresright]{rotating}

\journal{Frontiers of Mathematics in China (2016, To appear)}

\begin{document}

\begin{frontmatter}

\title{Domain of attraction of quasi-stationary distribution for one-dimensional diffusions}

\author[]{Hanjun Zhang}
\author[]{Guoman He\corref{cor1}}
\ead{hgm0164@163.com}
\cortext[cor1]{Corresponding author}
\address{School of Mathematics and Computational Science, Xiangtan University, Hunan 411105, PR China}

\begin{abstract}
We study quasi-stationarity for one-dimensional diffusions killed at 0, when 0 is a regular boundary and $+\infty$ is an entrance boundary. We give a necessary and sufficient condition for the existence of exactly one quasi-stationary distribution, and we also show that this distribution attracts all initial distributions.
\end{abstract}


\end{frontmatter}

\section{Introduction}
\label{sect1}
In this paper, we consider the one-dimensional diffusions $X$ on $[0,\infty)$ given by
\begin{equation}
\label{1.1}
dX_t=dB_t-q(X_t)dt,~~~~~~~ X_0=x>0,
\end{equation}
where $(B_t;t\geq0)$ is a standard one-dimensional Brownian motion and $q\in C^1([0,\infty))$. Observe that, under the condition $q\in C^1([0,\infty))$, $\int_{0}^{d}e^{Q(y)}dy<\infty$ and $\int_{0}^{d}e^{-Q(y)}dy<\infty$ for some (and, therefore, for all) $d>0$, which is equivalent to saying that the boundary point 0 is regular in the sense of Feller, where $Q(y)=\int_{0}^{y}2q(x)dx$.
\par
Let $\mathbb{P}_x$ and $\mathbb{E}_x$ stand for the probability and the expectation, respectively, associated with $X$ when initiated from $x$. For any distribution $\nu$ on $(0,\infty)$, we define $\mathbb{P}_\nu(\cdot):=\int_{0}^{\infty}\mathbb{P}_x(\cdot)\nu(dx)$. Let $\tau$ be the hitting time of 0, that is,
\begin{equation*}
\tau=\inf\{t>0: X_t=0\}.
\end{equation*}
Associated to $X$, we consider the sub-Markovian semigroup given by $T_tf(x)=\mathbb{E}_x(f(X_t),\tau>t)$, with density kernel denoted by $r(t,x,y)$. We denote by $L$ the infinitesimal operator of $X$, that is,
\begin{equation*}
L={1\over2}\partial_{xx}-q\partial_x.
\end{equation*}
Let us introduce the following useful measure defined on $(0,\infty)$:
\begin{equation}
\label{1.2}
\mu(dy):=e^{-Q(y)}dy.
\end{equation}
Notice that $\mu$ is the speed measure for $X$.
\par
For all the results in this paper, we will use the following hypothesis.
\vskip0.3cm
\noindent {\bf Hypothesis~(H).} $~~~~~~~~~\int_{0}^{\infty}e^{Q(y)}\left(\int_y^{\infty}e^{-Q(z)}dz\right)dy<\infty.$
\vskip0.3cm
\noindent It will be seen in Section 2 that hypothesis $(\mathrm{H})$ means that $+\infty$ is an entrance boundary in the sense of Feller.
\par
One of the fundamental problems for a killed Markov process conditioned on survival is to study its long-term asymptotic behavior. Conditional stationarity, which we call quasi-stationarity, is one of the most interesting topics in this direction. More formally, the following definition captures the main object of interest of this work.
\begin{defi}
\label{defi1}
We say that a probability measure $\nu$ supported  on $(0,\infty)$  is a quasi-stationary distribution $\mathrm{(QSD)}$, if for all $t\geq0$ and any Borel subset $A$ of $(0,\infty)$,
\begin{equation*}
\mathbb{P}_\nu(X_t\in A|\tau>t)=\nu(A).
\end{equation*}
\end{defi}
\par
Following \cite{CCLMMS09}, let us introduce a notion, which is closely related to QSD, so-called quasi-limiting distribution (QLD).
\begin{defi}
\label{defi2}
We say that a probability measure $\pi$ supported  on $(0,\infty)$ is a $\mathrm{QLD}$, if there exists a probability measure $\nu$ such that the following limit exists in distribution$:$
\begin{equation*}
\lim_{t\rightarrow\infty}\mathbb{P}_\nu(X_t\in\bullet|\tau>t)=\pi(\bullet).
\end{equation*}
We also say that $\nu$ is attracted to $\pi$, or is in the domain of attraction of $\pi$, for the conditional evolution.
\end{defi}
\par
In some cases, the long-time behavior of the conditional distribution can be proved that is independent
of the initial state. This leads us to study the notion of Yaglom limit.
\begin{defi}
\label{defi3}
We say that a probability measure $\pi$ supported  on $(0,\infty)$ is a Yaglom limit, if for any $x\in(0,\infty)$
\begin{equation*}
\lim_{t\rightarrow\infty}\mathbb{P}_x(X_t\in\bullet|\tau>t)=\pi(\bullet).
\end{equation*}
\end{defi}
\par
It is generally believed that QSD, QLD and Yaglom
limit have the following relation (see \cite{MV12}):
\begin{equation*}
\mathrm{Yaglom}~\mathrm{limit}\Longrightarrow \mathrm{QSD} \Longleftrightarrow \mathrm{QLD}.
\end{equation*}
\par
A complete treatment of the QSD problem for a given family of processes should accomplish the following two things (see \cite{PAG95}):
\vskip0.1cm
$(\mathrm{i})$ determination of all QSD's; and
\vskip0.1cm
$(\mathrm{ii})$ solve the domain of attraction problem, namely, characterize all laws $\upsilon$ such that a given QSD $\nu$ attracts all $\upsilon$.
\vskip0.1cm
Although ever since the pioneering work by Mandl \cite{M61}, the existence of the Yaglom limit and that of
a QSD for killed one-dimensional diffusion processes have been proved by many authors (see, e.g., \cite{CMS95, KS12, MM94, SE07}), it is very difficult to give a complete answer to the question of domains of attraction for initial distributions are different from the Dirac measures and the compactly supported initial distributions. In fact, details about $(\mathrm{ii})$ are known only for the Brownian motion with strictly negative constant drift \cite{MPM98} and the Ornstein-Uhlenbeck process \cite{LM00}. Under Mandl's conditions are not satisfied, the problem of the existence, uniqueness and domain of attraction of QSDs for one-dimensional diffusions killed at 0 and whose drift is allowed to go to $-\infty$ at 0 and the process is allowed to have an entrance boundary at $+\infty$, are solved in a satisfactory way by Cattiaux et~al. \cite{CCLMMS09}. In the present paper, we will show that there is exactly one QSD for  one-dimensional diffusions $X$ killed at 0, when 0 is a regular boundary and $+\infty$ is an entrance boundary, and that this distribution attracts all initial distributions.
\par
In this paper, we give a necessary and sufficient condition for the existence of exactly one QSD in terms of $q$ (see hypothesis $(\mathrm{H})$). If the ground state $\eta_1$ (eigenfunction associated to the bottom of the spectrum $\lambda_1$) belongs to $\mathbb{L}^1(\mu)$, we show that this unique QSD $\nu_1$ can be written by
\begin{equation*}
d\nu_1=\frac{\eta_1d\mu}{\langle\eta_1,1\rangle_\mu},
\end{equation*}
where $\langle f,g\rangle_\mu:=\int_{0}^{\infty}f(u)g(u)\mu(du)$. In order to obtain $\eta_1\in\mathbb{L}^1(\mu)$, we show that the spectrum of $L$ is discrete (see Proposition \ref{prop2} below).
Our main results are Theorems \ref{thm 3.1} and \ref{thm 4.1} below. Moreover, we point out an interesting fact that $\eta_1$ is bounded (see Proposition \ref{prop 4.2} below).
\par
This paper is organized as follows. In Section 2, we study the spectrum of the operator ${L}$. We show in Section 3 that there exists a unique QSD for the process $X$. In the last section, we solve the problem of domain of attraction of this unique QSD.

\section{Spectrum of ${L}$}
\label{sect2}
Throughout this paper, we shall assume the process $X$ has a finite lifetime, i.e., for $x>0$,
\begin{equation*}
\mathbb{P}_x(\tau<\infty)=1.
\end{equation*}
This is very closely related to the following scale function of $X$:
\begin{equation}
\label{2.1}
\Lambda(x)=\int_{0}^{x}e^{Q(y)}dy~~~~\mathrm{for}~x\in (0,\infty).
\end{equation}
In fact, as can be seen from the definition of natural scale that $\Lambda(X_t)$ is a nonnegative local martingale, and so that $\mathbb{P}_x(\tau<\infty)=1$ for $x>0$ if and only if the scale function is infinite at $\infty$; that is, for $c>0$, $\int_{c}^{\infty}e^{Q(y)}dy=\infty$.
\par
If $\int_{0}^{\infty}e^{Q(y)}\left(\int_0^{y}e^{-Q(z)}dz\right)dy=\infty$ and hypothesis $(\mathrm{H})$ holds, then $+\infty$ is called an entrance boundary according to Feller's classification (see, e.g., \cite[Chapter 15]{KT81}).
Observe that $\int_{0}^{\infty}e^{Q(y)}\left(\int_0^{y}e^{-Q(z)}dz\right)dy=\infty$ if $\Lambda(\infty)=\infty$.
In fact, for any $c_2\in(0,\infty)$, we have
\begin{eqnarray*}
&&\int_{0}^{\infty}e^{Q(y)}\left(\int_0^{y}e^{-Q(z)}dz\right)dy\\
&=&\int_{0}^{c_2}e^{Q(y)}\left(\int_{0}^{y}e^{-Q(z)}dz\right)dy+\int_{c_2}^{\infty}e^{Q(y)}\left(\int_{0}^{y}e^{-Q(z)}dz\right)dy\\
&\geq&\int_{c_2}^{\infty}e^{Q(y)}\left(\int_{0}^{c_2}e^{-Q(z)}dz\right)dy\\
&=&\infty.
\end{eqnarray*}
However, if (H) is satisfied, then we can deduce that $\Lambda(\infty)=\infty$ (see Lemma \ref{Lem 2.1} below). Hence, hypothesis (H) directly implies that $+\infty$ is an entrance boundary.
\begin{lem}
\label{Lem 2.1}
Assume that $(\mathrm{H})$ holds. Then $\mu(0,\infty)<\infty$ and $\Lambda(\infty)=\infty$.
\end{lem}
\begin{pf}
For $x>0$, since $q\in C^1([0,\infty))$, then $\int_{0}^{x}e^{-Q(z)}dz<\infty$. In addition, (H) implies that $\int_x^{\infty}e^{-Q(z)}dz<\infty$. We thus prove that $\mu(0,\infty)<\infty$.
By using the Cauchy--Schwarz inequality, we get $x^2=\left(\int_{0}^{x}e^{Q(z)/2}e^{-Q(z)/2}dz\right)^2\leq\int_{0}^{x}e^{Q(z)}dz\int_{0}^{x}e^{-Q(z)}dz$, and therefore, (H) implies that $\Lambda(\infty)=\infty$.\qed
\end{pf}
\par
For the study of QSDs, we need to consider the positivity of the bottom of the spectrum $\lambda_1$ for one-dimensional diffusions. We have the following result.
\begin{lem}
\label{Lem 2.2}
Assume that $(\mathrm{H})$ holds. Then $\lambda_1>0$.
\end{lem}
\begin{pf}
For any $x>0$, when $0< y\leq x$, we have
\begin{eqnarray*}
\int_{0}^{x}e^{Q(y)}dy\int_x^{\infty}e^{-Q(z)}dz&=&\int_{0}^{x}e^{Q(y)}\int_x^{\infty}e^{-Q(z)}dzdy\\
                                                &\leq&\int_{0}^{x}e^{Q(y)}\int_y^{\infty}e^{-Q(z)}dzdy\\
                                                &\leq&\int_{0}^{\infty}e^{Q(y)}\int_y^{\infty}e^{-Q(z)}dzdy.
\end{eqnarray*}
It follows that
$$\delta:=\sup_{x>0}\int_{0}^{x}e^{Q(y)}dy\int_x^{\infty}2e^{-Q(y)}dy\leq\int_{0}^{\infty}e^{Q(y)}\int_y^{\infty}2e^{-Q(z)}dzdy.$$
If (H) is satisfied, then $\delta<\infty$. We know from \cite[Theorem 1.1]{MFC00} or \cite[Theorem 1]{PRG09} that $(4\delta)^{-1}\leq\lambda_1\leq(\delta)^{-1}$. From this estimate, the result follows.\qed
\end{pf}
\par
On the spectrum of ${L}$, we have the following result.
\begin{prop}
\label{prop2}
Assume that $(\mathrm{H})$ holds. Then we have\\
$~~~~~~~(\mathrm{i})$ the spectrum of $L$ is discrete$;$\\
$~~~~~~(\mathrm{ii})$
for any nonnegative $f, g\in \mathbb{L}^2(\mu),$
\begin{equation}
\label{2.2}
\lim_{t\rightarrow\infty}e^{\lambda_1t}\langle g,T_tf\rangle_\mu=\langle\eta_1,f\rangle_\mu\langle\eta_1,g\rangle_\mu.
\end{equation}
\end{prop}
\begin{pf}
(i) When (H) is satisfied, we can deduce that
\begin{equation*}
\lim_{n\rightarrow\infty}\mu([r,\infty))\int_{n}^{\infty}e^{Q(x)}dx=0.
\end{equation*}
Since
\begin{equation*}
\sup\limits_{r>n}\mu([r,\infty))\int_{n}^{r}e^{Q(x)}dx\leq\mu([r,\infty))\int_{n}^{\infty}e^{Q(x)}dx,
\end{equation*}
so we get
\begin{equation}
\label{2.3}
\lim\limits_{n\rightarrow\infty}\sup\limits_{r>n}\mu([r,\infty))\int_{n}^{r}e^{Q(x)}dx=0.
\end{equation}
We know from \cite[Theorem 1]{PRG09} that $(\ref{2.3})$ is equivalent to $\sigma_{ess}({L})=\emptyset$, i.e., the spectrum of $L$ is discrete, where $\sigma_{ess}({L})$ denotes the essential spectrum of ${L}$.
\par
(ii) This was proved in \cite[Theorem 3.1]{LJ12} for a drifted Brownian motion killed at 0, when 0 is an exit boundary and $+\infty$ is an entrance boundary, but the proof carries over without essential changes to our case. It is straightforward from the $\mathbb{L}^2$ version of the process.\qed
\end{pf}

\section{Existence and uniqueness of a quasi-stationary distribution}
\label{sect3}
In this section, we study the standard QSDs of a one-dimensional diffusion $X$ killed at 0, when 0 is a regular boundary and $+\infty$ is an entrance boundary, a typical problem for absorbing Markov processes.
According to \cite{FOT10}, the one-dimensional diffusion process $X$ is symmetric with respect to $\mu$ and satisfies the condition of \cite[Lemma 6.4.5]{FOT10}, thus we know that there exists a ground state $\eta_1$ of the Dirichlet form $(\mathcal {E},\mathscr{D}(\mathcal {E}))$ uniquely up to a sign and $\eta_1$ can be taken to be strictly positive on $(0,\infty)$.
\par
In \cite[Theorem 4.14]{KS12}, it is proved that there exists a unique QSD for $X$ if and only if the boundary point $+\infty$ is entrance in the sense of Feller, namely, $\int_{0}^{\infty}e^{Q(y)}\left(\int_0^{y}e^{-Q(z)}dz\right)dy=\infty$ and hypothesis $(\mathrm{H})$ holds. Actually, in Section 2, we see that hypothesis $(\mathrm{H})$ implies $\int_{0}^{\infty}e^{Q(y)}\left(\int_0^{y}e^{-Q(z)}dz\right)dy=\infty$.
The following theorem is the main result of this section.
\begin{thm}
\label{thm 3.1}
The following statements are equivalent$:$\\
$~~~~~~(\mathrm{i})$ $(\mathrm{H})$ holds$;$\\
$~~~~~(\mathrm{ii})$ there exists a unique quasi-stationary distribution
\begin{equation*}
d\nu_1=\frac{\eta_1d\mu}{\langle\eta_1,1\rangle_\mu}
\end{equation*}
for the process $X$.
\end{thm}
\begin{pf}
$(\mathrm{i})\Rightarrow(\mathrm{ii})$. First, notice that Proposition \ref{prop2} implies $\eta_1\in\mathbb{L}^2(\mu)$. By using Lemma \ref{Lem 2.1} and the Cauchy--Schwarz inequality, we obtain $\eta_1\in\mathbb{L}^1(\mu)$.
Thanks to the symmetry of the semigroup, for all $f\in \mathbb{L}^2(\mu)$, we have
\begin{equation}
\label{3.6}
\int T_tf\eta_1 d\mu=\int fT_t\eta_1 d\mu=e^{-\lambda_1 t}\int f\eta_1 d\mu.
\end{equation}
Equality $(\ref{3.6})$ can be extended to all bounded function $f$. In particular, we may use it with $f={\bf1}_A$ and with $f={\bf1}_{(0,\infty)}$, where ${\bf1}_A$ is the indicator function of $A$. Note that
$$\int T_t{\bf1}_A\eta_1 d\mu=\mathbb{P}_{\nu_1}(X_t\in A,\tau>t){\langle\eta_1,1\rangle_\mu}$$
and
$$\int T_t({\bf1}_{(0,\infty)})\eta_1 d\mu=\mathbb{P}_{\nu_1}(\tau>t){\langle\eta_1,1\rangle_\mu},$$
then
\begin{eqnarray*}
\mathbb{P}_{\nu_1}(X_t\in A|\tau>t)&=&{{\mathbb{P}_{\nu_1}(X_t\in A,\tau>t)}\over{{\mathbb{P}_{\nu_1}(\tau>t)}}}={{\int T_t{\bf1}_A\eta_1 d\mu}\over{\int T_t({\bf1}_{(0,\infty)})\eta_1 d\mu}}\\
&=&{{\int {\bf1}_AT_t\eta_1d\mu}\over{\int ({\bf1}_{(0,\infty)})T_t\eta_1 d\mu}}={{e^{-\lambda_1 t}\int_A\eta_1 d\mu}\over{e^{-\lambda_1 t}\int_{0}^{\infty}\eta_1 d\mu}}\\
&=&\nu_1(A).
\end{eqnarray*}
Thus, we get that $\nu_1$ is a QSD. Moreover, we know from Theorem \ref{thm 4.1} below that $\nu_1$ is the unique QSD of $X$.
\par
$(\mathrm{ii})\Rightarrow(\mathrm{i})$. If there is a unique QSD for $X$, then it is easy to prove that $\int_{0}^{\infty}e^{Q(y)}\left(\int_y^{\infty}e^{-Q(z)}dz\right)dy<\infty$. In fact, if $\int_{0}^{\infty}e^{Q(y)}\left(\int_y^{\infty}e^{-Q(z)}dz\right)dy=\infty$, then we know from \cite[Theorem 3.1]{ZH16} that there is a one-parameter family of quasi-stationary distributions, which is a contradiction with the uniqueness. Thus, the result follows. \qed
\end{pf}

\section{Domain of attraction of quasi-stationary distribution}
\label{sect4}
In this section, we consider the problem of the domains of attraction in our framework. We have proved that $\nu_1$ is a QSD in previous section. Next, we will use the same arguments as in the proof of \cite[Theorem 5.2]{CCLMMS09} to show $\nu_1$ is the Yaglom limit distribution. To achieve this, we need to verify $r(1,x,y)\in\mathbb{L}^2(\mu)$. In fact, under the assumption that $(\mathrm{H})$ holds, for all $x>0$, there exists a locally bounded function $\theta(x)$ such that, for all $y>0$ and all $t\geq1$,
\begin{equation}
\label{4.1111}
r(t,x,y)\leq\theta(x)e^{-\lambda_1t}\eta_1(y).
\end{equation}
The proof of this property is similar to that of \cite[Lemma 5.3]{CCLMMS09} or \cite[Lemma 4.4]{KS12}, since it only uses $\eta_1\in\mathbb{L}^1(\mu)$ and we know from the proof of Theorem \ref{thm 3.1} that this fact is obtained easily. Thus, from (\ref{4.1111}), it is trivial to see that for $t=1$,
\begin{equation*}
\int_{0}^{\infty}r^2(1,x,y)\mu(dy)\leq \theta^2(x)e^{-2\lambda_1}\int_{0}^{\infty}\eta^2_1(y)\mu(dy)<\infty,
\end{equation*}
that is, $r(1,x,y)\in\mathbb{L}^2(\mu)$.
\begin{prop}
Assume that $(\mathrm{H})$ holds. Then for any $x>0$ and any Borel subset $A$ of $(0,\infty)$,
\begin{equation}
\label{3.4}
\lim\limits_{t\rightarrow\infty}e^{\lambda_1t}\mathbb{P}_x(\tau>t)=\eta_1(x)\langle\eta_1,1\rangle_\mu,
\end{equation}
\begin{equation}
\label{3.5}
\lim\limits_{t\rightarrow\infty}e^{\lambda_1t}\mathbb{P}_x(X_t\in A,\tau>t)=\nu_1(A)\eta_1(x)\langle\eta_1,1\rangle_\mu.
\end{equation}
This implies that
$$\lim\limits_{t\rightarrow\infty}\mathbb{P}_x(X_t\in A|\tau>t)=\nu_1(A),$$
that is, $\nu_1$ is the Yaglom limit distribution.
\end{prop}
\begin{pf}
When $(\mathrm H)$ is satisfied, we know from Lemma \ref{Lem 2.1} that $\mu$ is a bounded measure. For any Borel set $A\subseteq(0,\infty)$ such that ${\bf1}_A\in\mathbb{L}^2(\mu)$ and any $x>0$, $t>1$, we have
\begin{eqnarray*}
\mathbb{P}_x(X_t\in A,\tau>t)&=&\int\mathbb{P}_y(X_{t-1}\in A,\tau>t-1)r(1,x,y)\mu(dy)\\
                               &=&\int T_{t-1}({\bf1}_A)(y)r(1,x,y)\mu(dy)\\
                               &=&\int{\bf1}_A(y)(T_{t-1}r(1,x,\cdot))(y)\mu(dy).
\end{eqnarray*}
Since both ${\bf1}_A$ and $r(1,x,\cdot)$ are in $\mathbb{L}^2(\mu)$, by using Proposition \ref{prop2}, we obtain
\begin{equation*}
\lim\limits_{t\rightarrow\infty}e^{\lambda_1(t-1)}\mathbb{P}_x(X_t\in A,\tau>t)=\langle{\bf1}_A,\eta_1\rangle_\mu\langle r(1,x,\cdot),\eta_1\rangle_\mu.
\end{equation*}
Since
$$\int r(1,x,y)\eta_1(y)\mu(dy)=(T_1\eta_1)(x)=e^{-\lambda_1}\eta_1(x),$$
we get that $\nu_1$ is the Yaglom limit.\qed
\end{pf}
\par
In the next result, we give a sharper estimate on $\eta_1$, which has not been mentioned by previous authors (see, e.g., \cite{CCLMMS09, KS12, LJ12, SE07}). We will use the same arguments as in the proof of \cite[equality 7.3]{CCLMMS09} to show $e^{\lambda_1t}\mathbb{P}_x(\tau>t)$ is uniformly bounded in the variables $t$ and $x$.
\begin{prop}
\label{prop 4.2}
Assume that $(\mathrm{H})$ holds. Then $\eta_1$ is bounded.
\end{prop}
\begin{pf}
Let us first remark that for $0<x\leq x_0$, $\mathbb{P}_{x}(\tau>t)\leq\mathbb{P}_{x_0}(\tau>t)$. Thus, from (\ref{3.4}), we get that $\eta_1(x)\leq\eta_1(x_0)$.
On the other hand, if $(\mathrm H)$ is satisfied, based on \cite[Theorem 4.14]{KS12}, we can deduce that there is $x_0>0$ such that $B_1:=\sup_{x\geq x_0}\mathbb{E}_x[e^{\lambda_1\tau_{x_0}}]<\infty$. From (\ref{3.4}) again, we get that $B_2:=\sup_{u\geq 0}e^{\lambda_1u}\mathbb{P}_{x_0}(\tau>u)<\infty$. Then for $x>x_0$, we have
\begin{eqnarray*}
\mathbb{P}_{x}(\tau>t)&=&\int_{0}^{t}\mathbb{P}_{x_0}(\tau>u)\mathbb{P}_{x}(\tau_{x_0}\in d(t-u))+\mathbb{P}_{x}(\tau_{x_0}>t)\\
&\leq&B_2\int_{0}^{t}e^{-\lambda_1u}\mathbb{P}_{x}(\tau_{x_0}\in d(t-u))+\mathbb{P}_{x}(\tau_{x_0}>t)\\
&\leq&B_2e^{-\lambda_1t}\mathbb{E}_x[e^{\lambda_1\tau_{x_0}}]+e^{-\lambda_1t}\mathbb{E}_x[e^{\lambda_1\tau_{x_0}}]\\
&\leq&e^{-\lambda_1t}B_1(B_2+1).
\end{eqnarray*}
Thus, we get that $e^{\lambda_1t}\mathbb{P}_x(\tau>t)$ is uniformly bounded in the variables $t$ and $x$. By using (\ref{3.4}), it is easily seen that for $x>x_0>0$, $\eta_1(x)\leq\frac{B_1(B_2+1)}{\langle\eta_1,1\rangle_\mu}$. Hence, for any $x>0$, there exists $x_0>0$ such that $\eta_1(x)\leq\max\{\eta_1(x_0),\frac{B_1(B_2+1)}{\langle\eta_1,1\rangle_\mu}\}$. This completes the proof.\qed
\end{pf}
\par
Although, in the literature, there are several articles which have studied the problem of the domains of attraction (see, e.g., \cite{KS12, SE07}), our result is of particular interest in the analysis of the domain of attraction of QSDs for one-dimensional diffusions because we use the definition of QLD to prove the domain of attraction but not impose any condition on the initial distribution. Inspired by the proof of \cite[Lemma 19]{VD09}, we have the following result.
\begin{thm}
\label{thm 4.1}
Assume that $(\mathrm{H})$ holds. Then $\nu_1$ attracts all initial distributions $\nu$ supported in $(0,\infty)$, that is, for any Borel set $A\subseteq(0,\infty)$
\begin{equation*}
\lim_{t\rightarrow\infty}\mathbb{P}_\nu(X_t\in A|\tau>t)=\nu_1(A).
\end{equation*}
In particular, $\nu_1$ is the unique quasi-stationary distribution.
\end{thm}
\begin{pf}
Let $\nu$ be any fixed probability distribution whose support is contained in $(0,\infty)$.
When (H) is satisfied, we know from the proof of Proposition \ref{prop 4.2} that $e^{\lambda_1t}\mathbb{P}_x(\tau>t)$ is uniformly bounded in the variables $t$ and $x$, and $\eta_1$ is bounded. Then, by the dominated convergence theorem, one can integrate with respect to $\nu$ under the limit in (\ref{3.4}):
$$\lim\limits_{t\rightarrow\infty}e^{\lambda_1t}\mathbb{P}_\nu(\tau>t)=\int_{0}^{\infty}\eta_1(x)\int_{0}^{\infty}\eta_1(y)\mu(dy)\nu(dx).$$
The same holds for (\ref{3.5}):
$$\lim\limits_{t\rightarrow\infty}e^{\lambda_1t}\mathbb{P}_\nu(X_t\in A,\tau>t)=\nu_1(A)\int_{0}^{\infty}\eta_1(x)\int_{0}^{\infty}\eta_1(y)\mu(dy)\nu(dx).$$
This implies that
$$\lim\limits_{t\rightarrow\infty}\mathbb{P}_\nu(X_t\in A|\tau>t)=\lim\limits_{t\rightarrow\infty}\frac{\mathbb{P}_\nu(X_t\in A,\tau>t)}{\mathbb{P}_\nu(\tau>t)}=\nu_1(A).$$
We complete the proof.\qed
\end{pf}

\section*{Acknowledgements}
The first author would like to thank Professor Servet Mart\'{\i}nez for his kind hospitality during a visit to the Centro de Modelamiento Matem\'{a}tico of Universidad de Chile, where part of this work was done. The work is supported by the National Natural Science Foundation of China (Grant No. 11371301) and Hunan Provincial Innovation Foundation For Postgraduate (Grant No. CX2015B203).












\end{document}